\documentstyle[12pt]{article}
\font\teneufm=eufm10
\font\seveneufm=eufm7
\font\fiveeufm=eufm5
\newfam\frakturfam

\textfont\frakturfam=\teneufm
\scriptfont\frakturfam=\seveneufm
\scriptscriptfont\frakturfam=\fiveeufm

\newtheorem{lm}{Lemma}[section]
\newtheorem{theor}{Theorem}[section]
\newtheorem{co}{Corollary}[section]

\def\bee{\begin{eqnarray}}
\def\bes{\begin{eqnarray*}}
\def\eee{\end{eqnarray}}
\def\ees{\end{eqnarray*}}

\def\a{\alpha}
\def\b{\beta}
\def\te{\theta}

\def\s{\sigma}

\def\Proof{{\sl Proof.}\ }

\title{The Anick automorphism of free associative algebras}

\begin{document}
\date{}
\maketitle

\begin{center}

{\bf U.\,U.\,Umirbaev}\\
Eurasian National University,\\
 Astana, 010008, Kazakhstan \\
e-mail: {\em umirbaev@yahoo.com}

\end{center}

\begin{abstract}
We prove that the well-known Anick automorphism (see \cite[p. 343]{Cohn}) 
of the free associative algebra $F<x,y,z>$ over an arbitrary field 
$F$ of characteristic $0$ is wild.
\end{abstract}



\section{Introduction}

\hspace*{\parindent}

Let $B_n=F<y_1,y_2,\ldots,y_n>$ be the free associative algebra with free generators  
 $y_1,y_2,\ldots,y_n$ over a field $F$, and let  $Aut\,B_n$ be the automorphism group of $B_n$.
  Let  $\phi = (f_1,f_2,\ldots,f_n)$ denote an automorphism $\phi$ of $B_n$ such that 
   $\phi(y_i)=f_i,\, 1\leq i\leq n$. An automorphism  
\bes
\s(i,\a,f) = (y_1,\ldots,y_{i-1}, \a y_i+f, y_{i+1},\ldots,y_n),
\ees
where $0\neq\a\in F,\ f\in F<y_1,\ldots,y_{i-1},y_{i+1},\ldots,y_n>$, 
 is called {\em elementary}. The subgroup $TA(B_n)$ of $Aut\,B_n$ generated by all 
  elementary automorphisms is called the {\em tame automorphism group}, 
 and the elements of this subgroup are called the {\em tame automorphisms} 
 of $B_n$. Nontame automorphisms of $B_n$ are called {\em wild}.

\smallskip

In 1968 P.Cohn \cite{Cohn2} proved that the automorphisms of free 
Lie algebras with a finite set of generators are tame. 
 It is well known \cite{Czer,Jung,Kulk,Makar} that the automorphisms 
of polynomial algebras and free associative algebras in two variables are tame. 
 It was recently proved in    
\cite{SU4,Umi25,SU3} that the well-known Nagata automorphism (see \cite{Nagata})
\bes
\sigma&=&(x+(x^2-yz)z,\,y+2(x^2-yz)x+(x^2-yz)^2z,\,z)
\ees
of the polynomial algebra $F[x,y,z]$ over a field $F$ of characteristic $0$ is wild. 

In the case of free associative algebras the question about the 
existence of wild automorphisms was open and was formulated 
in \cite[Problem 1.74]{Dnestr} as a problem of P.\,Cohn. 
The best known candidate to be nontame was the Anick automorphism (see \cite[p. 343]{Cohn})
\bes
\delta&=&(x+z(xz-zy),\,y+(xz-zy)z,\,z)
\ees
of the algebra $F<x,y,z>$. 

 It is well known \cite{Smith} that the Nagata and the Anick automorphisms are stably tame.

This paper is devoted  to prove that the Anick automorphism   $\delta$ of the free associative algebra  
$F<x,y,z>$ over a field $F$ of characteristic $0$ is wild.  
In fact, we prove that the Anick automorphism   
$\delta$ induces a wild automorphism of the free metabelian associative algebra  in three variables. 

The paper is organized as follows. 
In Section 2 we recall the definition of Fox derivatives of free associative 
 algebras and formulate the main technical result of the paper 
(Theorem \ref{t1}). It was proved by P.\,Cohn \cite{Cohn0} that the matrix  
\bes
\left(\begin{array}{cc}
\frac{\partial \delta(x)}{\partial x} & \frac{\partial \delta(y)}{\partial x}\\
 \frac{\partial \delta(x)}{\partial y}& \frac{\partial \delta(y)}{\partial y}\\
\end{array}\right) 
=
\left(\begin{array}{cc}
1+z'\otimes z& 1\otimes z^2\\
-(z'\otimes 1)^2& 1-z'\otimes z\\
\end{array}\right)
\ees
is not a product of elementary matrices over the polynomial algebra $F[z'\otimes 1, 1\otimes z]$. 
From this fact we deduce that the Anick automorphism $\delta$ is wild. 

In Section 3 we begin to study automorphisms of the 
 free metabelian algebra $C$ 
 with free generators $z_1,z_2,z_3$. 
A description of defining relations of the 
tame automorphism group $TA(A)$ of the polynomial algebra $A=F[x_1,x_2,x_3]$ was given in \cite{UUU1}.   
Using this result, we get a characterization 
of tame automorphisms of $C$ modulo the tame automorphisms of A (Theorem \ref{t4.1}).

 In Section 4 we study the Jacobian matrices 
of tame automorphisms of the free metabelian associative algebra $C$. 
Given properties of these matrices imply the proof  of Theorem \ref{t1}.

The results of this paper were first published in \cite{UUU}.

\section{Main results}

\hspace*{\parindent}

First we recall the 
definition of Fox derivatives of free associative algebras (see, for example \cite{Dicks}). 
 Let $D$ be an arbitrary associative algebra (with unit 1) over $F$. We denote by 
$D'$ the algebra anti-isomorphic to $D$ and the corresponding anti-isomorphism we denote  by 
$':D\rightarrow D'$. The algebra $U(D)=D'\otimes_FD$ is called {\em the universal 
multiplicative enveloping algebra of $D$}. Let $I_D$ be the kernel of the mapping 
\bes
\lambda : D'\otimes_FD \rightarrow D, \ \ \lambda(\sum f_i'\otimes g_i)= \sum f_i' g_i.
\ees 
Then the mapping  
\bes
\Delta : D \rightarrow I_D, \ \ \Delta(f)= f'\otimes 1 - 1\otimes f, \ \ f\in D, 
\ees
is called {\em the universal derivation of $D$} \cite{Dicks}. Note that 
\bes
\Delta(fg)=\Delta(f)(1\otimes g)+\Delta(g)(f'\otimes 1), \ \ f,g\in D. 
\ees

Let us consider the free associative algebra $B_n=F<y_1,y_2,\ldots,y_n>$. 
In this case $I_{B_n}$ is a free right $U(B_n)$-module \cite{Dicks}  with free generators 
\bes
\Delta(y_i)=y_i'\otimes 1 - 1\otimes y_i, \ \ 1\leq i \leq n. 
\ees 
Then the Fox derivatives $\frac{\partial f}{\partial y_i}$ of an arbitrary element $f\in B_n$ can be uniquely defined by the formula
\bes
\Delta(f)=\sum_{i=1}^{n} \Delta(y_i) \frac{\partial f}{\partial y_i}. 
\ees
Note that the Fox derivatives  $\frac{\partial f}{\partial y_i}$ also can be defined as a linear map  
\bes
\frac{\partial}{\partial y_i}: B_n \rightarrow U(B_n), \ \ 1\leq i\leq n, 
\ees
such that   
\bee\label{f1}
\frac{\partial y_j}{\partial y_i}=\delta_{ij}, \ \ 
\frac{\partial fg}{\partial y_i}=\frac{\partial f}{\partial y_i} (1\otimes g)+\frac{\partial g}{\partial y_i} (f'\otimes 1), 
\eee
where $\delta_{ij}$ is Kronecker symbol.

For every $f\in B_n$ we put  
\bes
\partial(f)=\left(\frac{\partial f}{\partial y_1},\frac{\partial f}{\partial y_2},\ldots,\frac{\partial f}{\partial y_n}\right)^t, 
\ees
where $t$ means transposition. If $\varphi$ is an endomorphism of $B_n$, then the matrix  
\bes
J(\varphi)=(\partial(\varphi(y_1)),\partial(\varphi(y_2)),\ldots,\partial(\varphi(y_n)))
\ees
is called the Jacobian matrix of $\varphi$. For arbitrary endomorphisms   $\varphi$ and $\psi$ we have  
\bee\label{f2}
J(\varphi\psi)=J(\varphi)\varphi(J(\psi)).
\eee
Consequently, $J(\varphi)$ is invertible over $U(B_n)$ if $\varphi$ is an automorphism. The converse statement is also true \cite{Dicks,Schof}.

Let us consider automorphisms of the algebra $B=F<y_1,y_2,y_3>$. 
For an arbitrary automorphism $\varphi$ of $B$ we put  
\bes
J_2(\varphi)=
\left(\begin{array}{cc}
\frac{\partial \varphi(y_1)}{\partial y_1} & \frac{\partial \varphi(y_2)}{\partial y_1}\\
\frac{\partial \varphi(y_1)}{\partial y_2} & \frac{\partial \varphi(y_2)}{\partial y_2}\\
\end{array}\right).
\ees

Consider the homomorphism 
\bes
\pi : B=F<y_1,y_2,y_3> \longrightarrow A=F[x_1,x_2,x_3], \ \ \pi(y_i)=x_i, \ \ 1\leq i\leq 3. 
\ees
Note that the kernel $R$ of the homomorphism $\pi$ is the commutator ideal of $B$. 
If $\te=(f_1,f_2,f_3)\in Aut\,B$, then obviously  $\pi(\te)=(\pi(f_1),\pi(f_2),\pi(f_3))\in Aut\,A$. 
We say that the automorphism  $\pi(\te)$ is induced by $\te$.  

If $\phi : D\rightarrow H$ is a homomorphism of 
associative algebras, then the homomorphism $U(D)\rightarrow U(H)$ defined as 
$f'\otimes g\longmapsto\phi(f)'\otimes \phi(g)$ we will also denote by $\phi$. 

We now consider the homomorphism  
\bes
\nu  : B\longrightarrow F[y_3], \ \ \nu(y_1)=\nu(y_2)=0, \, \nu(y_3)=y_3.
\ees
Note that $U(F[y_3])= F[y_3'\otimes 1, 1\otimes y_3]$ is a polynomial algebra in 2 variables. 

By $e_{ij}$ denote the standard matrix units.  
A matrix $E_{ij}(d)=E+de_{ij}$, where $i\neq j$ and $d\in D$, is called an {\em elementary matrix} over $D$. 
 By $E_2(D)$ we denote the subgroup of $SL_2(D)$ generated by all elementary matrices.

Our main technical result is the following 
\begin{theor}\label{t1}
Assume that $\varphi$ is a tame automorphism of the free associative algebra $B=F<y_1,y_2,y_3>$ 
over an arbitrary field $F$ of characteristic $0$ such that $\pi(\varphi)=id$. 
Then 
\bes
\nu(J_2(\varphi))\in E_2(U(F[y_3])).
\ees 
\end{theor}

We prove this theorem in Section 4. Now we give some essential corollaries. 

\begin{co}\label{c1}
The Anick automorphism  
\bes
\delta&=&(x+z(xz-zy),\,y+(xz-zy)z,\,z)
\ees
of the free associative algebra $F<x,y,z>$ over a field $F$ of characteristic $0$ is wild.  
\end{co}
\Proof
We have 
\bes
\nu(J_2(\delta))= J_2(\delta) 
=
\left(\begin{array}{cc}
1+z'\otimes z& 1\otimes z^2\\
-(z'\otimes 1)^2& 1-z'\otimes z\\
\end{array}\right). 
\ees
It is well known \cite{Cohn0} that $\nu(J_2(\delta))\notin E_2(F[z'\otimes 1, 1\otimes z])$. 
Put $\delta\s=\psi$, where  
\bes
\s=\s(1,1,-y) \s(2,1,-xz^2) \s(1,1,y). 
\ees
We can immediately check that  
\bes
\pi(\psi)=id, \ \ \nu(J_2(\psi))=\nu(J_2(\delta)) \nu(\delta(J_2(\s))), \ \ 
\nu(\delta(J_2(\s)))\in E_2(F[z'\otimes 1, 1\otimes z]). 
\ees
Consequently, $\nu(J_2(\psi))\notin E_2(F[z'\otimes 1, 1\otimes z])$. 
It follows from Theorem \ref{t1} that $\psi\notin TA(B)$. Therefore $\delta\notin TA(B)$, i.e., $\delta$ is wild.  
$\Box$

\smallskip

\begin{co}\label{c2}
Let $\varphi=(f,g,z)$  be an automorphism of the free associative algebra  \\
$F<x,y,z>$ over a field $F$ of characteristic $0$ such that $\deg_{x,y}f, \deg_{x,y}g \leq 1$. 
Then $\varphi$  is tame if and only if 
\bes
J_2(\varphi) \in GL_2(F) E_2(F[z'\otimes 1, 1 \otimes z]). 
\ees  
\end{co}
\Proof
Assume that the matrix $J_2(\varphi)$ satisfies the condition of the corollary. 
Without loss of generality we can suppose that 
\bes
J_2(\varphi) \in E_2(F[z'\otimes 1, 1 \otimes z]). 
\ees  
Since $\deg_{x,y}f, \deg_{x,y}g \leq 1$, to each elementary matrix of $E_2(F[z'\otimes 1, 1 \otimes z])$ 
corresponds an elementary transformation of the first two coordinates of $\varphi$. 
Applying elementary transformations corresponding to a decomposition of $J_2(\varphi)$, 
we can get $(x+a(z), y+b(z), z)$ from $\varphi$. Consequently, $\varphi$ is tame. 

Now, assume that $\varphi$ is tame. Then $\varphi$ induces a tame automorphism of $A$. 
Since $\deg_{x,y}f, \deg_{x,y}g \leq 1$, it is not difficult to find
 a tame automorphism $\psi=(f',g',z)$ such that 
$\pi(\varphi \psi)=id$,  $\deg_{x,y}f', \deg_{x,y}g '\leq 1$, and 
\bes
J_2(\psi)\in GL_2(F) E_2(F[z'\otimes 1, 1 \otimes z]). 
\ees  
We have also $J_2(\varphi\psi)=J_2(\varphi)\varphi(J_2(\psi))$. Using Theorem \ref{t1}, we obtain     
\bes
J_2(\varphi\psi)\in GL_2(F) E_2(F[z'\otimes 1, 1 \otimes z]).  
\ees
 Consequently, $J_2(\varphi)$ satisfies the condition of the corollary. 
 $\Box$

\section{A characterization of tame automorphisms}

\hspace*{\parindent}

Let us remind that the variety of associative algebras defined by the identity  
\bes
[x,y][z,t]=0 
\ees
is called {\em metabelian} \cite{Umi5}. 
 This important variety of algebras are playing 
the role of the metabelian groups in combinatorial group theory. In \cite{Umi5} there were 
introduced and investigated the Jacobian matrices 
of endomorphisms of free metabelian associative algebras.

Let us denote by $C_n$ the free metabelian associative algebra (with unit 1)  with the set of free generators 
  $Z=\{z_1,z_2,\ldots,z_n\}$ and by $A_n=F[x_1,x_2,\ldots , x_n]$ the polynomial algebra in the variables 
  $x_1,x_2,\ldots , x_n$. By $C$ denote the free metabelian associative algebra  
with free generators $z_1,z_2,z_3$. 

It was proved in \cite{Umi5} that every 
automorphism of the polynomial algebra $A_n$ can be lifted to an automorphism 
of the algebra $C_n$. Then $C$ has wild automorphisms since the algebra $A$ has them.  Note that the Anick automorphism   
$\delta$ induces a tame automorphism of $A$. Later we see that $\delta$ induces a wild automorphism of $C$.

 The tame automorphism group 
 $TA(C_n)$ of $C_n$ is generated by all elementary automorphisms 
\bee\label{f3.1} 
\s(i,\a,f) = (z_1,\ldots,z_{i-1}, \a z_i+f, z_{i+1},\ldots,z_n),
\eee
where $0\neq\a\in F,\ f\in F<z_1,\ldots,z_{i-1},z_{i+1},\ldots,z_n>$. 
  
 The following relations between elementary automorphisms are known \cite{UUU1}: 
\bee\label{f3.2} 
\s(i,\a,f) \s(i,\b,g) = \s(i,\a \b, \b f+g);  
\eee
\bee\label{f3.3} 
\s(i,\a,f)^{-1} \s(j,\b,g) \s(i,\a,f) = \s(j,\b,\s(i,\a,f)^{-1}(g)),  
\eee
where $i\neq j$ and $f\in <Z\setminus \{z_i,z_j\}>$; 
\bee\label{f3.4}
\s(i,\a,f)^{(ks)}=\s(j,\a,(ks)(f)),
\eee
where $1\leq k\neq s\leq n$, $(ks)=\s(s,-1,z_k) \s(k,1,-z_s) \s(s,1,z_k)$, and $z_j=(ks)(z_i)$. 

\smallskip 

However the relations (\ref{f3.2})--(\ref{f3.4}) may fail to be a full system of defining relations of 
the group  $TA(C_n)$. 
It was proved in \cite{UUU1} that the relations (\ref{f3.2})--(\ref{f3.4}) 
are defining relations of the tame automorphism group $TA(A)$ 
with respect to the generators (\ref{f3.1}). 

\smallskip

Consider the homomorphism 
\bes
\tau : C \longrightarrow A=F[x_1,x_2,x_3], \ \ \tau(z_i)=x_i, \ \ 1\leq i\leq 3. 
\ees
Note that the kernel $I=Ker\,\tau$ of the homomorphism $\tau$ is the commutator ideal of $C$. 
According to \cite{Umid1}, elements of the type 
\bes
z_{j_1}z_{j_2}\ldots z_{j_k}[z_{i}z_{i_1}]z_{i_2}\ldots z_{i_s}, 
\ees 
where $j_1\leq j_2\leq \ldots \leq j_k$ and $i>i_1\leq i_2\leq \ldots \leq i_s$, is a basis of $I$.  

If $\te=(f_1,f_2,f_3)\in Aut\,C$, then $\tau^*(\te)=(\tau(f_1),\tau(f_2),\tau(f_3))\in Aut\,A$. 
 It is easy to check that the map  
\bes
\tau^*: Aut\,C \longrightarrow Aut\,A
\ees
is a group homomorphism. We denote by $T$ the restriction of $\tau^*$ to $TA(C)$, i.e.,   
\bes
T: TA(C) \longrightarrow TA(A), \ \ T(\te)=\tau^*(\te), \ \ \te\in TA(C). 
\ees

Obviously, any automorphism of the form   
\bee\label{f4.1}
\s(i,\a,f),\,\,\,\a=1,\,f\in I,  
\eee
is contained in $Ker\,T$. 

\begin{theor}\label{t4.1}
Let  $F$ be a field of characteristic $0$. Then the kernel of the homomorphism  
\bes
T: TA(C) \longrightarrow TA(A),  
\ees
as a normal subgroup of the group $TA(C)$, is generated by all automorphisms of the form (\ref{f4.1}).
\end{theor}
\Proof
Let $N$ be a normal subgroup of  $TA(C)$ generated by all elements of the form (\ref{f4.1}). 
Since $N\subseteq Ker\,T$, the homomorphism  $T$ induces the group homomorphism 
\bes
\mu: G=TA(C)/N\longrightarrow TA(A). 
\ees

We show that $\mu$ is an isomorphism. Every element $f\in C$ has a  unique representation  in the form $f=f_0+f_1$, 
where $f_0$ is a linear combination of the words $z_1^iz_2^jz_3^k$ and $f_1\in I$. We have  
\bes
\s(i,\a,f)=\s(i,\a,f_0)\s(i,1,f_1). 
\ees
Since $\s(i,1,f_1)\in N$, the group $G$ is generated by all elements  of the form 
\bee\label{f4.2}
\s(i,\a,f_0)N, \ \ f\in C. 
\eee
Consequently, $\mu$ establishes a one to one correspondence between the set of generators 
(\ref{f4.2}) of the group $G$ and the set of generators (\ref{f3.1}) of the group $TA(A)$. 
Of course, the relations (\ref{f3.2})--(\ref{f3.4}) are also true for the elementary automorphisms 
of $C$. They induce the same relations for the elements (\ref{f4.2}). 
According to \cite{UUU1}, the relations (\ref{f3.2})--(\ref{f3.4}) are defining 
relations of the group $TA(A)$ with respect to the generators (\ref{f3.1}). Then  there exists a homomorphism 
\bes
\omega  : TA(A) \longrightarrow G 
\ees
which is inverse to $\mu$ on the set of generators. Consequently, $\mu$ is an isomorphism. 
 $\Box$

For every $\varphi\in Aut\,C$ we denote by $\overline{\varphi}$ the induced 
automorphism of $A$, i.e. $\overline{\varphi}=\tau^*(\varphi)\in Aut\,A$. 

\begin{co}\label{c4.1}
Let  $\varphi$ be a tame automorphism of the algebra $C$ such that  $\overline{\varphi}=id$. Then 
\bes
\varphi=\phi_1^{\sigma_1}\phi_2^{\sigma_2}\ldots \phi_k^{\sigma_k}, 
\ees
where $\phi_i$ are automorphisms of the form (\ref{f4.1}), $\sigma_i\in TA(C)$, 
and $a^b=b^{-1}ab$ is the group conjugation.
\end{co}

\section{Jacobian matrices of tame automorphisms}

\hspace*{\parindent}

Let us recall the definition of Fox derivatives for the free metabelian algebra $C$ \cite{Umi5}. 
Consider the homomorphism 
\bes
\varepsilon : B \longrightarrow C, \ \ \varepsilon(z_i)=y_i, \ \ 1\leq i\leq 3. 
\ees
Note that $Ker(\varepsilon)=R^2$, where $R$ is the commutator ideal of $B$. 

Now consider the mappings 
\bes
\pi \frac{\partial}{\partial y_i} : B \longrightarrow U(A), \ \ 1\leq i\leq 3. 
\ees 
It follows from (\ref{f1}) that $R^2\subseteq Ker(\pi \frac{\partial}{\partial y_i})$. 
Then for any element $f\in C$ Fox derivatives $\frac{\partial f}{\partial z_i}$ 
can be well defined as $\frac{\partial f}{\partial z_i}=\pi(\frac{\partial g}{\partial y_i})$, 
where $\varepsilon(g)=f$ and $1\leq i\leq 3$. Note that $\frac{\partial f}{\partial z_i}\in U(A)$ 
and $U(A)$ is a polynomial algebra in six variables. 

For any endomorphism $\varphi$ of the algebra $C$ the Jacobian matrix $J(\varphi)$ and the matrix $J_2(\varphi)$ 
can be defined as in Section 2. Instead of (\ref{f2}), in this case for arbitrary endomorphisms   
$\varphi$ and $\psi$ of the algebra $C$ we have  
\bee\label{f100}
J(\varphi\psi)=J(\varphi) \overline{\varphi}(J(\psi)).
\eee
The Jacobian matrix $J(\varphi)$ is invertible over $U(A)$ 
iff $\varphi$ is an automorphism \cite{Umi5}.

 The commutator ideal $I$ of the algebra $C$ has a right $U(C)$-module structure with 
respect to the action  
\bes
m\cdot(f'\otimes g)= fmg, \ \ m\in I, f,g\in C. 
\ees
Since $I^2=0$, this action induces a right 
$U(A)$-module structure on $I$. If $f\in I$, then  $f$ can be written in the form  
\bee\label{f4.4}
f=[z_1,z_2]\cdot f_{12}+[z_1,z_3]\cdot f_{13}+[z_2,z_3]\cdot f_{23}, \ \ f_{ij}\in U(A). 
\eee
Using (\ref{f1}), we can easily get    
\bee\label{f4.5}
\frac{\partial f}{\partial z_3}=\Delta(x_1)f_{13}+\Delta(x_2)f_{23}. 
\eee

We denote by $e_i$ the vector-row of length 3 whose $i$th component is $1$ and the others are $0$, $1\leq i\leq 3$. 
\begin{lm}\label{l4.1} 
Let $\varphi$ be an automorphism of the algebra $C$ of the type (\ref{f4.1}) and let $\psi\in Aut\,C$. Then   
\bes
J(\psi\varphi\psi^{-1})=E+u^t\cdot v, 
\ees
where $u=(u_1,u_2,u_3)$, $v=(v_1,v_2,v_3)$, $u_i,v_i\in U(A)$. 
\end{lm}
\Proof
Assume that $\varphi=\s(i,1,f)$, where $f\in I$. Then $J(\varphi)=E+\partial(f)\cdot e_i$. We have  
\bes
J(\psi\varphi\psi^{-1})=J(\psi)\overline{\psi}(J(\varphi) \overline{\varphi}(J(\psi)))=J(\psi)\overline{\psi}(J(\varphi)) \overline{\psi}(J(\psi))\\
=J(\psi)\overline{\psi}(J(\varphi)) J(\psi)^{-1}=E+J(\psi)\overline{\psi}(\partial(f))\cdot \overline{\psi}(e_i) J(\psi)^{-1}.
\ees
Put $u^t=J(\psi)\overline{\psi}(\partial(f))$, $v=\overline{\psi}(e_i) J(\psi)^{-1}$. $\Box$

\smallskip

For an arbitrary automorphism $\varphi$ of $C$ we put  
\bes
J_2(\varphi)=
\left(\begin{array}{cc}
\frac{\partial \varphi(z_1)}{\partial z_1} & \frac{\partial \varphi(z_2)}{\partial z_1}\\
\frac{\partial \varphi(z_1)}{\partial z_2} & \frac{\partial \varphi(z_2)}{\partial z_2}\\
\end{array}\right).
\ees  

We now consider the homomorphism  
\bes
\eta : A\longrightarrow F[x_3], \ \ \eta(x_1)=\eta(x_2)=0, \, \eta(x_3)=x_3.
\ees

\begin{lm}\label{l4.3} 
Let $\varphi\in Aut\,C$ and $\overline{\varphi}=id$. Then   
\bes
\eta(J(\varphi))=
\left(\begin{array}{cc}
\eta(J_2(\varphi)) & \ast \\
0 & 1\\
\end{array}\right) 
\ees
\end{lm}
\Proof
We have $\varphi=(z_1+f_1,z_2+f_2,z_3+f_3)$, where $f_i\in I$. 
It follows from  (\ref{f4.5}) that $\eta(\frac{\partial f_i}{\partial z_3})=0$. $\Box$

\bigskip

For every $0\neq\lambda\in U(A)$ by $U_\lambda(A)$ denote the localization 
of $U(A)$ with respect to $\{\lambda^i\}_{i\geq 0}$. 
Note that $U_\lambda(A)$ is the subalgebra of the 
quotient algebra $Q(U(A))$ generated by $U(A)\cup\{\lambda^{-1}\}$. 

By 
$diag(\a,\b)$ denote the diagonal square matrix with the diagonal entries $\a,\b$.
 Now we put $S_\lambda=E_2(U_\lambda(A))\cdot diag(\lambda,1)$.  It is easy to check that  
\bes
E_{12}(a)\cdot diag(\lambda,1)=diag(\lambda,1)\cdot E_{12}(\frac{a}{\lambda}), \ \ 
E_{21}(a)\cdot diag(\lambda,1)=diag(\lambda,1)\cdot E_{21}(\lambda a), \\
diag(\lambda,1)=E_{12}(1)\cdot E_{21}(-1)\cdot E_{12}(1)\cdot diag(1,\lambda)\cdot E_{12}(-1)\cdot E_{21}(1)\cdot E_{12}(-1). 
\ees
From here we get  
\bee\label{f4.6}
S_\lambda=E_2(U_\lambda(A))\cdot diag(1,\lambda)=diag(\lambda,1)\cdot E_2(U_\lambda(A))=diag(1,\lambda)\cdot E_2(U_\lambda(A)).
\eee

\begin{lm}\label{l4.2} 
Let $\varphi$ be an automorphism of $C$ of the form (\ref{f4.1}) and let $\psi\in T(C)$. Then   
\bes
\lambda=\det(J_2(\psi\varphi\psi^{-1}))\neq 0, \ \ \ \ J_2(\psi\varphi\psi^{-1})\in S_\lambda.  
\ees
\end{lm}
\Proof 
Since $\overline{\psi\varphi\psi^{-1}}=id$, it follows that $\det(J(\psi\varphi\psi^{-1}))=1$. 
By Lemma \ref{l4.3}, we have $\eta(\lambda)=1$, i.e.   $\lambda\neq 0$. 

 Assume that $\varphi=\s(i,1,f)$, where $f\in I$. Then 
\bes
J(\varphi)=E+\partial(f)\cdot e_i.
\ees

 If $i=3$, then $J_2(\varphi)=E$. Since $f$ does not contain $z_i$, the  $i$th component of $\partial(f)$ is $0$. 
 Consequently, if $i\leq 2$, then the matrix $J_2(\varphi)$ is elementary. 
 Since $\det(J_2(\varphi))=1$, the statement of the lemma is valid for $\varphi$.  

We have  
\bes
\s(i,\a,g)=(1i) \s(1,\a, (1i)(g)) (1i)  
\ees
if $i>1$. Since  
\bes
(1i)=\s(1,-1,0) \s(1,1,z_i) \s(i,1,-z_1) \s(1,1,z_i), 
\ees
the relations (\ref{f3.2}) and (\ref{f3.3}) give  
\bes
\s(i,\a,g)=\s(1,1,z_i) \s(i,1,z_1) \s(1,\a,(1-\a)z_i-(1i)(g)) \s(i,1,-z_1) \s(1,1,z_i). 
\ees
Therefore we can assume that $\psi$ is a product of the automorphisms   
\bee\label{f4.7}
\s(1,\a,g), \ \ X_{ij}(\b)=\s(j,1,\b x_i), \ \ 1<j, \ \ i\neq j. 
\eee
Consider $\psi=\psi_1\psi_2\ldots \psi_k$, where $\psi_i$, $1\leq i\leq k$, are automorphisms of the form (\ref{f4.7}). 
We now begin an induction on  $k$. Above we have checked that the lemma is true if $k=0$.  

Assume that the statement of the lemma is valid for 
\bes
\varphi_1=\psi_2\ldots \psi_k \varphi (\psi_2\ldots \psi_k)^{-1}, 
\ees
i.e. 
$J_2(\varphi_1)\in S_\lambda$, where $\lambda=\det(J_2(\varphi_1))$. By Lemma \ref{l4.2}, we have 
\bes
J(\varphi_1)=E+u^t\cdot v,
\ees
 where $u=(u_1,u_2,u_3)$, $v=(v_1,v_2,v_3)$. Then 
\bes
J_2(\varphi_1)=E+(u_1,u_2)^t\cdot (v_1,v_2)\in S_\lambda.
\ees 
 Consequently, $\overline{\psi}_1(J_2(\varphi_1))\in S_\mu$, where $\mu=\overline{\psi}_1(\lambda)=\det(\overline{\psi}_1(J_2(\varphi_1)))$. 

Put  
\bes
\overline{\psi}_1(u)=(b_1,b_2,b_3)=b, \ \ \overline{\psi}_1(v)=(c_1,c_2,c_3)=c,
\ees
\bes
P=\overline{\psi}_1(J_2(\varphi_1)), \ \ Q=J_2(\psi_1\varphi_1\psi_1^{-1}).
\ees
 Then  
\bes
P=E+(b_1,b_2)^t\cdot (c_1,c_2). 
\ees
Since $\overline{\varphi}_1=id$, it follows that  
\bes
J(\psi_1\varphi_1\psi_1^{-1})=J(\psi_1) \overline{\psi}_1(J(\varphi_1)) J(\psi_1)^{-1}=
E+J(\psi_1)b^t\cdot cJ(\psi_1)^{-1}.
\ees

We first consider the case when $\psi_1=\s(1,\a,g)$. We have  
\bes
\partial(\psi_1(z_1))=(\a,a_2,a_3)^t, \ \ J(\psi_1)=(\partial(\psi_1(z_1)),e_2^t,e_3^t). 
\ees
An immediate calculation gives  
\bes
Q=E+(\a b_1,a_2 b_1+b_2)^t\cdot (\a^{-1}(c_1-c_2 a_2),c_2)+(\a b_1,a_2 b_1+b_2)^t\cdot (-\a^{-1}c_3 a_3, 0). 
\ees
Hence   
\bes
Q_1=E+(\a b_1,a_2 b_1+b_2)^t\cdot (\a^{-1}(c_1-c_2 a_2),c_2)\\
=diag(\a,1)\cdot E_{21}(a_2)\cdot P\cdot E_{21}(-a_2)\cdot  diag(\a^{-1},1).  
\ees
 It follows from (\ref{f4.6}) that $Q_1\in S_\mu$. By (\ref{f4.6}), we can also assume that $Q_1=P_1D$, where $P_1\in E_2(U_\mu(A))$, $D=diag(1,\mu)$. Then  
\bes
Q=P_1(D+P_1^{-1}(\a b_1,a_2 b_1+b_2)^t\cdot (-\a^{-1}c_3 a_3, 0)). 
\ees
We have  
\bes
D+P_1^{-1}(\a b_1,a_2 b_1+b_2)^t\cdot (-\a^{-1}c_3 a_3, 0)
=
\left(\begin{array}{cc}
\lambda_1& 0\\
\lambda_2& \mu\\
\end{array}\right).  
\ees
It follows  from (\ref{f4.6}) that  
\bes
D+P_1^{-1}(\a b_1,a_2 b_1+b_2)^t\cdot (-\a^{-1}c_3 a_3, 0)\in S_\nu, 
\ees
where $\nu=\lambda_1 \mu =\det(Q)$. Since $\mu|\nu$, we have $U_\mu(A)\subseteq U_\nu(A)$ and $E_2(U_\mu(A))\subseteq E_2(U_\nu(A))$. 
Consequently, $Q\in S_\nu$. 

We now assume that $\psi_1=X_{ij}(\b)$, $j>1$, $i\neq j$. We have $J(\psi_1)=E_{ij}(\b)$. Suppose that $i=1$. If $j=2$, then  
\bes
Q=E+(b_1+\b b_2, b_2)^t\cdot 
(c_1, c_2-\b c_1)
=E_{12}(\b)\cdot P\cdot E_{12}(-\b)\in S_\mu. 
\ees
If $j=3$, then  
\bes
Q=E+(b_1+\b b_3, b_2)^t\cdot 
(c_1, c_2)
=P+(\b b_3, 0)^t\cdot 
(c_1, c_2). 
\ees
By (\ref{f4.6}), we can assume that $P=DP_1$, where $P_1\in E_2(U_\mu(A))$, $D=diag(1,\mu)$. Then  
\bes
Q=(D+(\b b_3, 0)^t\cdot 
(c_1, c_2)P_1^{-1})P_1
 \in S_\nu, \ \ \mu|\nu, \ \ \nu=\det(Q). 
\ees

Suppose that $i=2$. Then $j=3$ and  
\bes
Q=E+(b_1,b_2+\b b_3)^t\cdot 
(c_1, c_2)
=P+(0,\b b_3)^t\cdot 
(c_1, c_2). 
\ees
By (\ref{f4.6}), in this case we can assume that $P=DP_1$, where $P_1\in E_2(U_\mu(A))$, $D=diag(\mu,1)$. Hence  
\bes
Q=(D+(0,\b b_3)^t\cdot 
(c_1, c_2)P_1^{-1}) P_1\in S_\nu, \ \ \mu|\nu, \ \ \nu=\det(Q). 
\ees

Now suppose that $i=3$. Then $j=2$ and  
\bes
Q=E+(b_1, b_2)^t\cdot 
(c_1, c_2-\b c_3)
=P+(b_1, b_2)^t\cdot 
(0, -\b c_3). 
\ees
Put $P=P_1D$, where $P_1\in E_2(U_\mu(A))$, $D=diag(\mu,1)$. Then we have  
\bes
Q=P_1(D+P_1^{-1}(b_1, b_2)^t\cdot 
(0, -\b c_3))\in S_\nu, \ \ \mu|\nu, \ \ \nu=\det(Q). 
\ees
Thus $J_2(\psi\varphi\psi^{-1})=Q\in S_\nu$. $\Box$

\begin{lm}\label{l4.4} 
Let $\varphi$ be an automorphism of $C$ of the form (\ref{f4.1}) and let $\psi\in T(C)$. Then   
\bes
\eta(J_2(\psi\varphi\psi^{-1}))\in E_2(U(F[x_3])).  
\ees
\end{lm}
\Proof 
Assume that $\det(J_2(\psi\varphi\psi^{-1}))=\lambda$. By Lemma \ref{l4.3}, we have $\eta(\lambda)=1$. 
Consequently, the homomorphism $\eta$  extends to a homomorphism of the algebra 
$U_\lambda(A)$ to $U(F[x_3])$.  Moreover, we have $\eta(S_\lambda)\subseteq E_2(U(F[x_3]))$. 
Now Lemma \ref{l4.2} gives our statement.  
$\Box$

\bigskip

\begin{theor}\label{t4.2}
Assume that $\varphi$ is a tame automorphism of $C$ such that $\overline{\varphi}=id$. Then 
\bes
\eta(J_2(\varphi))\in E_2(U(F[x_3])).
\ees 
\end{theor}
\Proof
According to Theorem \ref{t4.1}, we have  
\bes
\varphi=\te_1\te_2\ldots \te_s, \ \ \te_i=\psi_i\varphi_i\psi_i^{-1}, 
\ees
where $\varphi_i$ is an automorphism of the form (\ref{f4.1}) and $\psi_i$ is a tame automorphism. 
Since $\overline{\te}_i=id$, applying (\ref{f100}) we have  
\bes
J(\varphi)=J(\te_1) J(\te_2)\ldots J(\te_s).  
\ees
Lemma \ref{l4.3} implies that  
\bes
\eta(J_2(\varphi))=\eta(J_2(\te_1)) \eta(J_2(\te_2))\ldots \eta(J_2(\te_s)).  
\ees
By Lemma \ref{l4.4}, we have $\eta(J_2(\te_i))\in E_2(U(F[x_3]))$.
Hence $\eta(J_2(\varphi))\in E_2(U(F[x_3]))$. 
 $\Box$

\bigskip

Using this theorem and the same discussions as in the proof of Corollary \ref{c1}, we get 
\begin{co}\label{c4.1}
The Anick automorphism $\delta$ induces a wild automorphism 
of the algebra $C$ over a field $F$ of characteristic $0$.  
\end{co}

\bigskip

{\bf Proof of Theorem \ref{t1}.} 
 Since  $\varphi=(f_1,f_2,f_3)$ is a tame automorphism of the algebra $B=F<y_1,y_2,y_3>$ 
and $\pi(\varphi)=id$, it follows that 
\bes
\psi =\varepsilon^*(\varphi)=(\varepsilon(f_1),\varepsilon(f_2),\varepsilon(f_3))\in TA(C)
\ees
 and $\overline{\psi}=id$. By Theorem \ref{t4.2} we have 
\bes
\nu(J_2(\psi))\in E_2(U(F[x_3])). 
\ees

Consider the isomorphism 
\bes
\rho : F[y_3]\longrightarrow F[x_3], \ \ \rho(y_3)=x_3. 
\ees
Then $\rho \nu =\eta \varepsilon$. Consequently, 
\bes
(\rho\nu)(J_2(\varphi)) =\rho(\nu(J_2(\varphi))) =\eta(\varepsilon(J_2(\varphi)))=\eta(J_2(\psi))\in E_2(U(F[x_3])). 
\ees
We have proved that  
\bes
\nu(J_2(\varphi))\in E_2(U(F[y_3]))
\ees
since $\rho$ is an isomorphism. $\Box$

\bigskip

\begin{center}
{\bf\large Acknowledgments} 
\end{center}

\hspace*{\parindent}

I am grateful to I.\,Shestakov for thoroughly going over the details
of the proofs. I am also grateful to Max-Planck Institute f\"ur Mathematik for hospitality and exellent working conditions. I also thank J.\,Alev, P.\,Cohn, N.\,Dairbekov, and L.\,Makar-Limanov for numerous helpful comments and discussions. 

\bigskip

\hspace*{\parindent}

\end{document}